\documentclass[a4paper,10pt]{scrartcl}
\usepackage[english]{babel}
\usepackage[utf8]{inputenc}
\usepackage{enumerate}
\usepackage{geometry}
\usepackage{amsfonts,amsmath,amssymb,amsopn}
\usepackage{amsthm}
\usepackage{mathtools}
\usepackage{stmaryrd}
\usepackage{nicefrac}
\usepackage{exscale}

\usepackage[numbers,square]{natbib}
\usepackage{url} 
\usepackage[colorlinks=true,pdfpagelabels,unicode]{hyperref}
\hypersetup{urlcolor=red, citecolor=blue} 

\allowdisplaybreaks 

\theoremstyle{plain}

\newtheoremstyle{theo}
	{3pt} 
	{3pt} 
	{\itshape} 
	{} 
		{\bfseries} 
	{\\} 
	{ } 
	{\thmname{#1}\thmnumber{ #2.}\thmnote{ - #3}} 
\theoremstyle{theo}

\newtheorem{definition}{Definition}[section]
\newtheorem{lemma}[definition]{Lemma}
\newtheorem{theorem}[definition]{Theorem}

\newenvironment{bew}{\begin{proof}[\bfseries Proof:]}{\end{proof}}

\DeclareMathOperator{\bomega}{\overline{\Omega}}
\DeclareMathOperator{\romega}{\partial\Omega}

\DeclareMathOperator{\supp}{supp}
\DeclareMathOperator{\intd}{d\!}

\DeclareMathOperator{\dive}{\nabla\cdot}
\DeclareMathOperator{\wto}{\rightharpoonup}

\newcommand{\epsi}{\varepsilon}

\usepackage[usenames,dvipsnames,svgnames,table]{xcolor}

\newcommand{\Tme}{T_{max,\;\!\epsi}}

\newcommand{\into}[1]{\int_0^{#1}\!}
\newcommand{\intot}{\into{t}}
\newcommand{\intoT}{\into{T}}
\newcommand{\intoinf}{\into{\infty}}

\newcommand{\intomega}{\int_{\Omega}\!} 

\newcommand{\intoTomega}{\intoT\!\intomega}
\newcommand{\intinfomega}{\int_0^\infty\!\!\intomega}

\newcommand{\intromega}{\int_{\romega}\!} 

\newcommand{\Lo}[1][1]{L^{#1}(\Omega)} 
\newcommand{\W}[1][1,2]{W^{#1}(\Omega)}

\newcommand{\LSp}[2]{L^{#1\;\!}\!\left(#2\right)} 
\newcommand{\LSploc}[2]{L_{loc}^{#1}\!\left(#2\right)} 

\newcommand{\CSp}[2]{C^{#1}\!\left(#2\right)}

\newcommand{\R}{\mathbb{R}}
\newcommand{\N}{\mathbb{N}}

\author{Tobias Black\thanks{Institut f\"ur Mathematik, Universit\"at Paderborn, Warburger Str. 100, 33098 Paderborn, Germany; email: \mbox{tblack@math.upb.de}}}
\title{Global generalized solutions to a parabolic-elliptic Keller-Segel system with singular sensitivity}

\begin{document}
\maketitle
\begin{abstract}
\noindent
{\textbf{Abstract:} We investigate the parabolic-elliptic Keller-Segel model
\begin{align*}
\left\{
\begin{array}{r@{\,}l@{\quad}l@{\quad}l@{\,}c}
u_{t}&=\Delta u-\,\chi\nabla\!\cdot(\frac{u}{v}\nabla v),\ &x\in\Omega,& t>0,\\
0&=\Delta v-\,v+u,\ &x\in\Omega,& t>0,\\
\frac{\partial u}{\partial\nu}&=\frac{\partial v}{\partial\nu}=0,\ &x\in\romega,& t>0,\\
u(&x,0)=u_0(x),\ &x\in\Omega,&
\end{array}\right.
\end{align*}
in a bounded domain $\Omega\subset\mathbb{R}^n$ $(n\geq2)$ with smooth boundary.

\noindent We introduce a notion of generalized solvability which is consistent with the classical solution concept, and we show that whenever $0<\chi<\frac{n}{n-2}$ and the initial data satisfy only certain requirements on regularity and on positivity, one can find at least one global generalized solution.
}\\[0.1cm]


\end{abstract}


\newpage
\section{Introduction}\label{sec1:intro}
Since the introduction of the original parabolic-parabolic Keller--Segel model (\cite{KS70})
\begin{align*}
\left\{
\begin{array}{r@{\,}l@{\quad}l@{\quad}l@{\,}c}
u_{t}=\Delta u&-\,\chi\nabla\!\cdot(u\nabla v),&&\\
v_{t}=\Delta v&-\,v+u,&&
\end{array}\right.
\end{align*}
cross-diffusive systems of this type have been prototypical models for the description of the biological phenomenon of chemotaxis, a process of self-enhanced migration of cells towards higher concentration of a signal substance. For an overview of the biological background and related models in the context of chemotaxis we refer the reader to the surveys \cite{HP09} and \cite{Ho03}.

In this work we will consider a parabolic-elliptic Keller--Segel system with logarithmic sensitivity, as described by
\begin{align}\label{KSsing}
\left\{
\begin{array}{r@{\,}l@{\quad}l@{\quad}l@{\,}c}
u_{t}&=\Delta u-\,\chi\nabla\!\cdot(\frac{u}{v}\nabla v),\ &x\in\Omega,& t>0,\\
0&=\Delta v-\,v+u,\ &x\in\Omega,& t>0,\\
\frac{\partial u}{\partial\nu}&=\frac{\partial v}{\partial\nu}=0,\ &x\in\romega,& t>0,\\
u(&x,0)=u_0(x),\ &x\in\Omega,&
\end{array}\right.
\end{align}
where $\Omega\subset\R^n$, $n\geq2$, is a bounded domain with smooth boundary and $\chi$ is a positive parameter. The singular sensitivity governing the cross-diffusive motion in the form featured in \eqref{KSsing} expresses the model assumption that stimulus perception is governed by the Weber--Fechner law (\cite{HP09},\cite{ROSEN1978}). Both the parabolic-elliptic system and the parabolic-parabolic variant 
\begin{align}\label{parapara}
\left\{
\begin{array}{r@{\,}l@{\quad}l@{\quad}l@{\,}c}
u_{t}=\Delta u&-\,\chi\nabla\!\cdot(\frac{u}{v}\nabla v),&&\\
\tau v_{t}=\Delta v&-\,v+u,&&
\end{array}\right.
\end{align}
have been intensively studied in the last decades and a large amount of literature has been dedicated to the investigation of conditions ensuring the existence of time-global solutions in contrast to the possible occurrence of solutions blowing up in finite time. Nevertheless, for general dimensions the permissible strength of the sensitivity, as measured by the parameter $\chi$, which allows for any solution stemming from reasonably regular initial data to exist globally, still remains mostly unclear even in the simple parabolic-elliptic setting of \eqref{KSsing}.

Let us briefly summarize some known results. In \cite{NagSen_JAMA98} Nagai and Senba studied radially symmetric solutions to \eqref{KSsing} and showed that the classical radially symmetric solutions are global and bounded, whenever the conditions $\chi>0$ and $n=2$, or $\chi<\frac{n}{n-2}$ and $n\geq3$ are fulfilled. On the other hand, if $n\geq3$ and $\chi>\frac{2n}{n-2}$ they could prove the existence of solutions blowing up in finite time. In the studies independently undertaken in \cite{Biler_JAMA99} under the condition $\chi<\frac{2}{n}$ the existence of global weak solutions was verified without any symmetry requirements, though the boundedness of the solutions was left open. Several years later this result was extended by proving that \eqref{KSsing} possesses unique global bounded solutions if $\chi<\frac{2}{n}$ (\cite{FWY15}) and instead of $\chi(v)=\frac{\chi}{v}$ even more general sensitivity functions satisfying $\chi(v)\leq\frac{\chi_0}{v^k}$ have been investigated (cf. also \cite{FWY15}). More recently, the question concerning bounded classical solutions has been solved for the case of $n=2$ and the existence of finite time blow-up has been completely ruled out for any $\chi>0$ (\cite{FujieSenba16}). In the corresponding parabolic-parabolic setting of \eqref{parapara} for suitably small $\tau\in(0,1)$, the results for two-dimensional domains do not differ significantly and still blow-up does not occur for any $\chi>0$ (cf. \cite{FujSen_Non16}). On the contrary when $\tau=1$, the possible choices for $\chi$ are slightly more restricted. For instance, in \cite{Winkler_MMAS11} the global existence of classical solutions was established for $\chi<\sqrt{\frac{2}{n}}$ and the boundedness of these solutions was later proven in \cite{Fujie_JMAA15} and recently generalized to sensitivity functions of the form $\chi(v)\leq\frac{\chi}{(a+v)^k}$ with $a\geq0$, $k\geq1$ and $\chi<k(a+\eta)^{k-1}\sqrt{\frac{2}{n}}$, for some $\eta>0$ possibly depending on the initial data when $k>1$, by the studies in \cite{MizuYoko2017}. That $\chi=\sqrt{\frac{2}{n}}$ is not the critical value in the parabolic-parabolic version of \eqref{KSsing} is illustrated by the results of \cite{Lan_MMAS16}, where for $n=2$ global classical solutions were obtained for $\chi<\chi_0$ with some $\chi_0>1$.

By weakening the solution concept, larger ranges for $\chi$ allowing for global solutions could be achieved in the setting of \eqref{parapara}. Global weak solutions are known to exist for $\chi<\sqrt{\frac{n+2}{3n-4}}$ (\cite{Winkler_MMAS11}) and in even weaker concepts global generalized solutions exist for $\chi<\sqrt{\frac{n}{n-2}}$ (\cite{StinnerWinkler_NARWA11}) or, as the most recent studies show, $\chi<\frac{n}{n-2}$, atleast when $n=2$ or $n\geq4$ (\cite{LanWin2017}), highlighting once more the importance of the case $\chi=\frac{n}{n-2}$ as witnessed in the radial parabolic-elliptic setting of \eqref{KSsing} in \cite{NagSen_JAMA98}. Since the result of \cite{LanWin2017} does not rely on any symmetry assumptions or largeness assumptions on the initial data, it seems reasonable to expect that in the simpler setting of \eqref{KSsing} global generalized solutions may also exist for $\chi<\frac{n}{n-2}$ without requiring either radial symmetric initial data or $n=2$.
\\[0.1cm]
\noindent{\textbf{Main results.} Our main purpose in this work is to introduce a concept of generalized solvability for \eqref{KSsing}, which on one hand is consistent with the concept of classical solvability and on the other hand is weak enough to allow for the construction of global solutions for $\chi<\frac{n}{n-2}$, without requiring any symmetry assumption on the initial data. To be precise, we will only assume the initial data $u_0$ to satisfy 
\begin{align}\label{IR}
u_0\in\CSp{0}{\bomega},\ u_0\geq0\ \text{in}\ \Omega,\quad \text{and}\quad\intomega u_0>0.
\end{align}
Under this assumption our main result can be stated as follows.
\begin{theorem}\label{theo:globgensol}
Let $n\geq 2$ and $\Omega\subset\R^n$ be a bounded domain with smooth boundary. Suppose that $$0<\chi<\frac{n}{n-2},$$ then for any $u_0$ satisfying \eqref{IR}, the problem \eqref{KSsing} possesses at least one global generalized solution $(u,v)$ in the sense of Definition \ref{def:glob-gen-sol} below. Furthermore, this solution satisfies
\begin{align}\label{eq:mass}
\intomega u(\cdot,t)=\intomega u_0\quad\text{for a.e. }t>0.
\end{align}
\end{theorem}

\setcounter{equation}{0} 
\section{Generalized solution concept}\label{sec2:gensol}
Concepts of generalized solvability in related settings have previously been studied in e.g. \cite{win15_chemorot}, \cite{LanWin2017}, \cite{Win16CS1}, or \cite{TB_stokes-sing2017}. The requirements imposed on the solution components in these concepts can often be viewed as generalizations of a classical supersolution properties combined with suitable regularity conditions giving meaning to the integral inequality prescribed. In the current setting the generalized solutions we will investigate will fulfill such a supersolution property for the quantity $u^p$, whereas for $v$ we will require the standard weak solution concept with everything made precise by the following definition.
\begin{definition}\label{def:glob-gen-sol}
Suppose that 
\begin{align}\label{eq:sol_def_reg1}
u\in\LSploc{1}{\bomega\times[0,\infty)}\qquad\text{and}\qquad v\in\LSploc{1}{[0,\infty);\W[1,1]}
\end{align}
are such that $u>0$ and $v>0$ a.e. in $\Omega\times(0,\infty)$. Then $(u,v)$ will be called a global generalized solution of \eqref{KSsing} if there exists $p\in(0,1)$ such that
\begin{align}\label{eq:sol_def_reg2}
\arraycolsep=1.4pt\def\arraystretch{1.25}
\left\{\begin{array}{l@{\, }c@{\, }l}
\nabla u^\frac{p}{2}&\in&\LSploc{2}{\bomega\times[0,\infty)},\\
u^{p+1}v^{-1}&\in&\LSploc{1}{\bomega\times[0,\infty)},\\
u^\frac{p}{2} v^{-1}\nabla v&\in&\LSploc{2}{\bomega\times[0,\infty)},
\end{array}\right.
\end{align}
and such that the inequality
\begin{align}\label{eq:sol_def_uineq}
-\intinfomega u^p&\varphi_t-\intomega u_0^p\varphi(\cdot,0)\nonumber\\
&\geq\frac{4(1-p)(1-p\chi)}{p}\intinfomega|\nabla u^\frac{p}{2}|^2\varphi+\intinfomega u^p\Delta\varphi+4(1-p)\chi\intinfomega\Big|\nabla u^\frac{p}{2}-\frac{u^\frac{p}{2}}{2v}\nabla v\Big|^2\varphi\nonumber\\
& \ \;-(1-p)\chi\intinfomega u^p\varphi+(1-p)\chi\intinfomega\frac{u^{p+1}}{v}\varphi-(1-2p)\chi\intinfomega u^p\frac{\nabla v}{v}\cdot\nabla\varphi
\end{align}
holds for each nonnegative $\varphi\in C^\infty_0\left(\bomega\times[0,\infty)\right)$ such that $\frac{\partial\varphi}{\partial \nu}=0$ on $\romega\times(0,\infty)$, if moreover the identity
\begin{align}\label{eq:sol_def_veq}
-\intinfomega\nabla v\cdot\nabla \psi+\intinfomega v\psi=\intinfomega u\psi
\end{align}
is valid for any $\psi\in C_0^\infty\!\left(\bomega\times[0,\infty)\right)$, and if $u$ satisfies 
\begin{align}\label{eq:sol_def_mass_ineq}
\intomega u(x,t)\intd x\leq\intomega u_0(x)\intd x\quad\text{for a.e. }t>0
\end{align}
as well as
\begin{align*}
u^\frac{p}{2}>0\qquad\qquad\text{a.e. in }\romega\times(0,\infty).
\end{align*}
\end{definition}

Let us now make sure that the above concept of solvability is consistent with the classical one, meaning that a generalized solution of \eqref{KSsing}, which is smooth enough, is also a classical solution of \eqref{KSsing}. In particular the a.e. positivity of $u^\frac{p}{2}$ on $\romega\times(0,\infty)$ will play a crucial role in making sure that $u$ satisfies its equation in the classical sense. The proof builds on ideas previously used in \cite[Lemma 2.1]{win15_chemorot} and \cite[Lemma 2.5]{LanWin2017}. 
\begin{lemma}\label{lem:compability-class-sol}
Let $\chi>0$ and assume that $(u,v)\in\Big(\CSp{0}{\bomega\times[0,\infty)}\cap\CSp{2,1}{\bomega\times(0,\infty)}\Big)\times\CSp{2,0}{\bomega\times(0,\infty)}$ is a global generalized solution of \eqref{KSsing} in the sense of Definition \ref{def:glob-gen-sol}. Then $(u,v)$ solves \eqref{KSsing} in the classical sense in $\Omega\times(0,\infty)$.
\end{lemma}

\begin{bew}
In light of the assumed regularity properties of $v$ it can be easily verified by standard arguments that $v$ is a classical solution to the second equation in \eqref{KSsing} with the prescribed initial and boundary data and we may focus on proving that $u$ is a classical solution of the first equation in \eqref{KSsing} for the remainder of the proof. Given an arbitrary nonnegative $\psi\in \CSp{\infty}{\bomega}$ satisfying $\frac{\partial\psi}{\partial\nu}\big|_{\romega}=0$, for $\epsi\in(0,1)$ we define $\varphi(x,t):=\psi(x)(1-\frac{t}{\epsi})_+$, $(x,t)\in\Omega\times(0,\infty)$ and plug $\varphi$ into \eqref{eq:sol_def_uineq} to see upon taking $\epsi\searrow0$ that
\begin{align*}
\intomega u^p(\cdot,0)\psi-\intomega u_0^p\psi\geq0\quad\text{for all nonnegative}\ \psi\in\CSp{\infty}{\bomega}\text{ with }\frac{\partial\psi}{\partial\nu}\big|_{\romega}=0,
\end{align*}
in view of Lebesgue's dominated convergence theorem and the continuity of $t\mapsto\intomega u^p(\cdot,t)$ at $t=0$. This readily establishes $u(\cdot,0)\geq u_0$ in $\Omega$, which in combination with the continuity of $u$ at $t=0$ and \eqref{eq:sol_def_mass_ineq} shows $u(\cdot,0)=u_0$ in $\Omega$.

Due to terms of the form $u^{p-1}$ appearing when we integrate by parts we will now only consider test functions $\varphi$ which are compactly supported in $\{u>0\}:=\{(x,t)\in\bomega\times[0,\infty)\,|\, u(x,t)>0\}$. Consequently, for any nonnegative $\varphi\in C_0^\infty(\bomega\times(0,\infty))$ with $\supp\varphi\subset\{u>0\}$ and $\frac{\partial\varphi}{\partial\nu}\big|_{\romega}=0$ integrating by parts in the second and last integrals on the right in \eqref{eq:sol_def_uineq} shows that
\begin{align}\label{eq:compability-proof1}
-&\intinfomega u^p\varphi_t-\intinfomega u_0^p\varphi(\cdot,0)\nonumber\\
&\ \ \geq \frac{4(1-p)(1-p\chi)}{p}\intinfomega|\nabla u^\frac{p}{2}|^2\varphi+\intinfomega \Big(p u^{p-1}\Delta u -\frac{4(1-p)}{p}|\nabla u^\frac{p}{2}|^2\Big)\varphi+p\intoinf\!\intromega u^{p-1}\frac{\partial u}{\partial\nu}\varphi
\nonumber\\&\quad\ +4(1-p)\chi\intinfomega\Big|\nabla u^\frac{p}{2}-\frac{u^\frac{p}{2}}{2v}\nabla v\Big|^2\varphi-(1-p)\chi\intinfomega u^p\varphi+(1-p)\chi\intinfomega\frac{u^{p+1}}{v}\varphi\\&\quad\ +(1-2p)\chi\intinfomega\Big(p \frac{u^{p-1}}{v}\nabla u\cdot\nabla v+u^p\frac{\Delta v}{v}-u^p\frac{|\nabla v|^2}{v^2}\Big)\varphi,\nonumber
\end{align}
since $\frac{\partial v}{\partial\nu}=0$ on $\romega\times(0,\infty)$. Now, straightforward calculations show that
\begin{align*}
\ &\frac{4(1-p)(1-p\chi)}{p}|\nabla u^\frac{p}{2}|^2+4(1-p)\chi\Big|\nabla u^\frac{p}{2}-\frac{u^\frac{p}{2}}{2v}\nabla v\Big|^2+(1-2p)p\chi\frac{u^{p-1}}{v}\nabla u\cdot\nabla v\\
=\ &\frac{4(1-p)}{p}|\nabla u^\frac{p}{2}|^2-2(1-p)p\chi\frac{u^{p-1}}{v}\nabla u\cdot\nabla v+(1-p)\chi\frac{u^p}{v^2}|\nabla v|^2+(1-2p)p\chi\frac{u^{p-1}}{v}\nabla u\cdot\nabla v\\=\ &
\frac{4(1-p)}{p}|\nabla u^\frac{p}{2}|^2+(1-p)\chi\frac{u^p}{v^2}|\nabla v|^2-p\chi\frac{u^{p-1}}{v}\nabla u\cdot\nabla v
\end{align*}
and thus, we may rewrite \eqref{eq:compability-proof1} as
\begin{align}\label{eq:compability-proof2}
\intinfomega \partial_t(u^p)\varphi&\geq\intinfomega p u^{p-1}\Big(\Delta u-\chi\frac{\nabla u\cdot\nabla v}{v}-\chi\frac{u}{v}\Delta v+\chi u\frac{|\nabla v|^2}{v^2}\Big)\varphi+p\intoinf\!\intromega u^{p-1}\frac{\partial u}{\partial\nu}\varphi\nonumber\\
&\quad+(1-p)\chi\intinfomega\frac{u^p}{v}\big(\Delta v-v+u\big)\varphi\nonumber\\
&=\intinfomega p u^{p-1}\Big(\Delta u-\chi\dive\Big(\frac{u}{v}\nabla v\Big)\Big)\varphi+p\intoinf\!\intromega u^{p-1}\frac{\partial u}{\partial\nu}\varphi
\end{align}
for any nonnegative $\varphi\in C_0^\infty(\bomega\times(0,\infty))$ with $\supp\varphi\subset\{u>0\}$ and $\frac{\partial\varphi}{\partial\nu}\big|_{\romega}=0$, since we already know that $v$ solves the second equation in \eqref{KSsing}. Restricting to nonnegative $\varphi\in C_0^\infty(\Omega\times(0,\infty)\cap\{u>0\})$ we may rely on a Du Bois-Reymond lemma type argument to conclude that
\begin{align}\label{eq:compability-proof3}
u_t\geq \Delta u-\chi\dive\Big(\frac{u}{v}\nabla v\Big)\qquad\text{in }\{u>0\}.
\end{align}
In view of the continuity of $u$ and the fact that $u>0$ a.e. in $\Omega\times(0,\infty)$, \eqref{eq:compability-proof3} actually holds in all of $\Omega\times(0,\infty)$.

In order to first see that $\frac{\partial u}{\partial\nu}\geq0$ on $\romega\times(0,\infty)$ we refer the reader to the proof of \cite[Lemma 2.5]{LanWin2017} for a detailed construction of suitable permissible test functions in \eqref{eq:compability-proof2} and remark here only that it is essential to ensure that the support of the test functions intersects the boundary only in points where $u$ is positive. Lastly, to show that $u$ fulfills the homogeneous Neumann boundary condition and that \eqref{eq:compability-proof3} is actually an equality, we integrate \eqref{eq:compability-proof3} over $\Omega\times(0,t)$ and make use of \eqref{eq:sol_def_mass_ineq} to see that 
\begin{align*}
\intomega u_0\geq\intomega u(\cdot,t)\geq\intomega u_0+\intot\!\intomega\Delta u-\chi\intot\!\intomega\dive\Big(\frac{u}{v}\nabla v\Big)=\intomega u_0+\intot\!\intromega\frac{\partial u}{\partial\nu}
\end{align*}
in view of Gauss' theorem and the fact that $\frac{\partial v}{\partial\nu}$ on $\romega\times(0,\infty)$. This shows that actually $\frac{\partial u}{\partial\nu}=0$ on $\romega\times(0,\infty)$ and in turn proves that \eqref{eq:compability-proof3} is an equality, implying that $u$ solves the first equation of \eqref{KSsing} in the classical sense.
\end{bew}

\setcounter{equation}{0} 
\section{Approximate solutions and basic properties}\label{sec3:approx}
The construction of a global generalized solution is based on a limit procedure of solutions to suitably regularized problems. We will therefore continue by investigating approximate problems which for $\epsi\in(0,1)$ take the form
\begin{align}\label{KSeps}
\left\{
\begin{array}{r@{\,}l@{\quad}l@{\quad}l@{\,}c}
u_{\epsi t}&=\Delta u_\epsi-\,\chi\nabla\!\cdot\big(\frac{u_\epsi}{(1+\epsi u_\epsi)v_\epsi}\nabla v_\epsi\big),\ &x\in\Omega,& t>0,\\
0&=\Delta v_\epsi-\,v_\epsi+u_\epsi,\ &x\in\Omega,& t>0,\\
\frac{\partial u_\epsi}{\partial\nu}&=\frac{\partial v_\epsi}{\partial\nu}=0,\ &x\in\romega,& t>0,\\
u_\epsi(&\!x,0)=u_0(x),\ &x\in\Omega.&
\end{array}\right.
\end{align}

\subsection{Local existence and first estimates independent of \texorpdfstring{$\epsi$}{epsilon}}\label{sec31:locex}
\begin{lemma}\label{lem:locex-approx}
Let $\chi>0$, $\epsi\in(0,1)$ and suppose that $u_0$ satisfies \eqref{IR}. Then there exists a maximal existence time $\Tme\in(0,\infty]$ and a unique pair $(u_\epsi,v_\epsi)$ of nonnegative functions 
\begin{align*}
u_\epsi&\in\CSp{0}{\bomega\times[0,\Tme)}\cap\CSp{2,1}{\bomega\times(0,\Tme)}\\
v_\epsi&\in\CSp{2,0}{\bomega\times(0,\Tme)}
\end{align*}
solving the problem \eqref{KSeps} in the classical sense. Moreover,
\begin{align}\label{eq:locex_alt}
\text{either }\quad\Tme=\infty\qquad\text{or}\qquad\|u_\epsi(\cdot,t)\|_{\Lo[\infty]}\to\infty\quad\text{as }\ t\nearrow\Tme.
\end{align}
\end{lemma}
\begin{bew}
The local existence of classical solutions on $\bomega\times(0,\Tme)$ and an extensibility criterion can be proven by relying on well-known fixed point arguments as displayed for a very closely related setting in \cite[Proposition 3.1]{FWY15} (or \cite[Lemma 3.1]{BBWT15} for a parabolic-parabolic variant), while making use of the facts that $\frac{u_\epsi}{1+\epsi u_\epsi}\leq\frac{1}{\epsi}$ and that $v_\epsi$ is strictly positive on $\bomega\times(0,\infty)$ (see Lemma \ref{lem:lowerbound-veps} below).
\end{bew}

In the sequel of the paper we will always assume that the initial data $u_0$ satisfy \eqref{IR}, and for $\epsi\in(0,1)$ we let $(u_\epsi,v_\epsi)$ denote the corresponding solution to \eqref{KSeps} given by Lemma \ref{lem:locex-approx}. We will first focus our efforts on proving that these local solutions are in fact global solutions to \eqref{KSeps}. As a starting point for further a priori estimates we obtain the following $L^1$--\,regularity result for the approximate solutions.
\begin{lemma}\label{lem:l1-bounds}
Let $\chi>0$ and $\epsi\in(0,1)$. Then
\begin{align*}
\intomega u_\epsi(\cdot,t)=\intomega u_0\quad
\qquad \text{and}\qquad
\intomega v_\epsi(\cdot,t)=\intomega u_0\qquad\text{for all }t\in(0,\Tme).
\end{align*}
\end{lemma}

\begin{bew}
The first asserted identity immediately follows from integration of the first equation in \eqref{KSeps}. Making use of the established mass conservation, an integration of the second equation in \eqref{KSeps} consequently proves the second equality.
\end{bew}

Another important property of the solutions to the approximate problems is a pointwise lower bound -- strictly larger than zero -- for the component $v_\epsi$, which can be shown by an estimation of the fundamental solution from below. Arguments of this in this spirit have previously been employed in e.g. \cite[Lemma 2.1]{FWY15} and \cite[Lemma 2.1]{MizuYoko2017}.
\begin{lemma}\label{lem:lowerbound-veps}
There exists $K_1>0$ such that for each $\epsi\in(0,1)$,
\begin{align*}
v_\epsi(x,t)\geq K_1\quad\text{for all }x\in\Omega\text{ and }t\in(0,\Tme).
\end{align*}
\end{lemma}

\begin{bew}
Making use of the positivity of the fundamental solution of the heat equation (e.g. \cite[Chapter 10]{Ito_diffusion_equations}), one can find $C_1>0$ such that for all nonnegative $\psi\in\CSp{0}{\bomega}$ we have
\begin{align*}
e^{t\Delta}\psi(x)\geq C_1\intomega\psi>0\quad\text{for all }x\in\Omega\text{ and all }t\in(0,\Tme).
\end{align*}
Now, we can make use of the representation of resolvents via semigroups and Lemma \ref{lem:l1-bounds} to obtain
\begin{align*}
v(x,t)=\intoinf e^{-t}e^{t\Delta}u(x,t)\intd t\geq C_1\intomega u_0\intoinf e^{-t}\intd t=C_1\intomega u_0=:K_1>0
\end{align*}
for all $x\in\Omega$ and all $t\in(0,\Tme)$.
\end{bew}

Additionally, we can make use of standard elliptic theory to slightly improve our a priori knowledge on the regularity of $v_\epsi$.
\begin{lemma}\label{lem:veps-bounds}
Let $q,r\geq1$ be such that $q<\frac{n}{n-2}$ and $r<\frac{n}{n-1}$. Then there exists $C>0$ such that for all $\epsi\in(0,1)$ we have
\begin{align}\label{eq:veps-bounds-veps}
\intomega v_\epsi^q(\cdot,t)\leq C\quad\text{for all }t\in(0,\Tme)
\end{align}
and
\begin{align}\label{eq:veps-bounds-nab}
\intomega|\nabla v_\epsi(\cdot,t)|^r\leq C\quad\text{for all }t\in(0,\Tme).
\end{align}
\end{lemma}

\begin{bew}
According to known results concerning elliptic boundary-value problems with inhomogeneities in $\Lo[1]$ (see e.g. \cite{BrezisStrauss}) one can find $C_1>0$ such that
\begin{align*}
\|v_\epsi(\cdot,t)\|_{\W[1,r]}\leq C_1\|-\Delta v_\epsi(\cdot,t)+v_\epsi(\cdot,t)\|_{\Lo[1]}+C_1\|v_\epsi(\cdot,t)\|_{\Lo[1]}\quad\text{for all }t\in(0,\Tme).
\end{align*}
In view of the second equation and Lemma \ref{lem:l1-bounds} this yields
\begin{align*}
\|v_\epsi(\cdot,t)\|_{\W[1,r]}\leq C_1\|u_\epsi(\cdot,t)\|_{\Lo[1]}+C_1\|v_\epsi(\cdot,t)\|_{\Lo[1]}\leq 2C_1\intomega u_0\quad\text{for all }t\in(0,\Tme),
\end{align*}
proving \eqref{eq:veps-bounds-nab}. This bound at hand, \eqref{eq:veps-bounds-veps} follows from the Sobolev embedding theorem, since for all $q\in[1,\frac{n}{n-2})$ we can pick some $r'\in[1,\frac{n}{n-1})$ such that $1-\frac{n}{r'}\geq-\frac{n}{q}$.
\end{bew}

In addition to the regularity provided by the previous result, we can also make rely on the strict positivity of $v_\epsi$ to obtain the following a priori bound in a straightforward manner.
\begin{lemma}\label{lem:ln-nab-veps_l2}
Let $\epsi\in(0,1)$. Then
\begin{align*}
\intomega\frac{|\nabla v_\epsi(\cdot,t)|^2}{v_\epsi(\cdot,t)^2}\leq |\Omega|\quad\text{for all }t\in(0,\Tme).
\end{align*}
\end{lemma}
\begin{bew}
In view of the strict positivity of $v_\epsi$ established by Lemma \ref{lem:lowerbound-veps}, we may use $v_\epsi^{-1}$ as a test function in the second equation of \eqref{KSeps} and obtain upon integration by parts that
\begin{align*}
\intomega\frac{|\nabla v_\epsi|^2}{v_\epsi^2}-\intomega\frac{v_\epsi}{v_\epsi}+\intomega\frac{u_\epsi}{v_\epsi}=0\quad\text{in }(0,\Tme).
\end{align*}
Due to the positivity of both $u_\epsi$ and $v_\epsi$, the assertion immediately follows upon dropping the integral containing $u_\epsi$.
\end{bew}

\subsection{Global solvability of the approximate problems}\label{sec32:globex}

Relying on the lower bound of $v_\epsi$, the fact that $\frac{u_\epsi}{1+\epsi u_\epsi}\leq\frac{1}{\epsi}$ holds for all $\epsi\in(0,1),$ and standard elliptic theory, we will make use of an iterative argument to improve the regularity of $u_\epsi$ and $\nabla v_\epsi$ to a level where semigroup arguments for the heat semigroup become applicable to provide the boundedness of $\|u_\epsi\|_{\Lo[\infty]}$, which in view of the extensibility criterion is sufficient to conclude that $\Tme=\infty$. 
\begin{lemma}\label{lem:globex-approx}
Let $\chi>0$ and $\epsi\in(0,1)$, and let $(u_\epsi,v_\epsi)$ denote the local classical solution of \eqref{KSeps} in $\Omega\times(0,\Tme)$ obtained in Lemma \ref{lem:locex-approx}. Then $\Tme=\infty$.
\end{lemma}

\begin{bew}
Suppose $\Tme\leq T_\epsi<\infty$. For $q>2$, by using integration by parts and Young's inequality we compute
\begin{align*}
\frac{1}{q}\frac{\intd}{\intd t}\intomega u_\epsi^q&=-(q-1)\intomega u_\epsi^{q-2}|\nabla u_\epsi|^2+(q-1)\chi\intomega\frac{u_\epsi^{q-1}}{(1+\epsi u_\epsi)v_\epsi}\nabla u_\epsi\cdot\nabla v_\epsi\\
&\leq C_1\intomega\frac{u_\epsi^q}{(1+\epsi u_\epsi)^2v_\epsi^2}|\nabla v_\epsi|^2\quad\text{for all }t\in(0,\Tme),
\end{align*}
with $C_1=\frac{(q-1)\chi^2}{4}>0$. Employing Young's inequality once more and making use of Lemma \ref{lem:lowerbound-veps} and the fact that $\frac{u_\epsi}{1+\epsi u_\epsi}\leq\frac{1}{\epsi}$ for all $\epsi\in(0,1)$ we obtain that
\begin{align}\label{eq:globapproxexeq1}
\frac{1}{q}\frac{\intd}{\intd t}\intomega u_\epsi^q\leq \intomega u_\epsi^q+\frac{C_2}{\epsi^q K_1^q}\intomega|\nabla v_\epsi|^q\quad\text{for all }t\in(0,\Tme),
\end{align}
with $C_2=\frac{(q-1)^2\chi^4}{16}>0$ and with $K_1>0$ given by Lemma \ref{lem:lowerbound-veps}. Furthermore, by standard elliptic theory (e.g. \cite[Theorem 9.32]{Brezis_funcana}) there exists $C_3>0$ such that 
\begin{align}\label{eq:globapproxexeq2}
\|v_\epsi(\cdot,t)\|_{\W[2,q]}\leq C_3\|u_\epsi(\cdot,t)\|_{\Lo[q]}\quad\text{holds for all }t>0,
\end{align}
since $q>2$, and combination with \eqref{eq:globapproxexeq1} shows that
\begin{align*}
\frac{\intd}{\intd t}\intomega u_\epsi^q\leq C_4\intomega u_\epsi^q,\quad\text{for all }t\in(0,\Tme)
\end{align*}
with $C_4=q\big(1+\frac{C_2C_3}{\epsi^q K_1^q}\big)>0$, implying that for any $q\in[1,\infty)$ we have $\|u_\epsi\|_{\Lo[q]}\leq C_5$ for all $t\in(0,\Tme)$ with some $C_5=C_5(\epsi,T_\epsi)>0$. Relying on \eqref{eq:globapproxexeq2} once more, we also see that $\|v_\epsi\|_{\W[1,\infty]}\leq C_6$ with some $C_6=C_6(\epsi,T_\epsi)>0$, due to the Sobolev embedding theorem. Thus, making use of a Moser type iteration (e.g. \cite[Lemma A.1]{TaoWin-quasilinear_JDE12}) we obtain $\|u_\epsi\|_{\Lo[\infty]}\leq C_7$ for all $t\in(0,\Tme)$ with some $C_7=C_7(\epsi,T_\epsi)>0$, contradicting the extensibility criterion \eqref{eq:locex_alt} and thereby proving that $\Tme=\infty$.
\end{bew}

\setcounter{equation}{0} 
\section{Construction of limit functions}\label{sec4:limit}
In the next section we will derive a fundamental inequality for the approximate systems \eqref{KSeps}. Relying on the fairly arbitrary choices possible for the test functions $\varphi$ used therein, we will then first apply this to $\varphi\equiv1$ (see Lemma \ref{lem:st-bounds}) to derive a set of crucial a priori estimates. Later on (see proof of Theorem \ref{theo:globgensol}) we will make use of this inequality to verify the supersolution property featured in \eqref{eq:sol_def_uineq} of Definition \ref{def:glob-gen-sol}.

\subsection{Precompactness properties}\label{sec4-1:precompact}
Before we start with the derivation of the fundamental inequality, let us introduce the following notation. For $p\in(0,1)$ and $\epsi\in(0,1)$ we define
\begin{align}\label{eq:Phi-def}
\Phi_\epsi(s):=p\int_0^s\frac{\sigma^{p-1}}{1+\epsi \sigma}\intd \sigma,\quad s\geq0.
\end{align}
Let us also remark here that obviously $\Phi_\epsi(s)\leq s^p$ for all $s\geq0$ and that for $\epsi\searrow0$ we have $\Phi_\epsi(s)\to s^p$, as these are two properties we will require later on.

\begin{lemma}\label{lem:testing-equation}
Let $\epsi\in(0,1)$, $\chi>0$, $p\in(0,1)$ and $T>0$. Assume that $\varphi\in\CSp{\infty}{\bomega\times[0,T]}$ satisfies $\frac{\partial\varphi}{\partial\nu}=0$ on $\romega\times(0,T)$. Then the classical solution $(u_\epsi,v_\epsi)$ of \eqref{KSeps} in $\Omega\times(0,\infty)$ satisfies
\begin{align}\label{eq:testing-eq}
-&\intoTomega u_\epsi^p\varphi_t+\intomega u_\epsi^p(\cdot,T)\varphi(\cdot,T)-\intomega u_\epsi^p(\cdot,0)\varphi(\cdot,0)\nonumber\\
&\ \geq\frac{4(1-p)(1-p\chi)}{p}\intoTomega|\nabla u_\epsi^\frac{p}{2}|^2\varphi+\intoTomega u_\epsi^p\Delta\varphi+4(1-p)\chi\intoTomega\Big|\nabla u_\epsi^\frac{p}{2}-\frac{u_\epsi^\frac{p}{2}}{2v_\epsi}\nabla v_\epsi\Big|^2\varphi\nonumber\\&\ \ \;
+(1-p)\chi\intoTomega\Phi_\epsi(u_\epsi)\varphi-2(1-p)\chi\intoTomega u_\epsi^p\varphi+(1-p)\chi\intoTomega\frac{u_\epsi^{p+1}}{v_\epsi}\varphi\\
&\ \ \ +(1-p)\chi\intoTomega\frac{\Phi_\epsi(u_\epsi)}{v_\epsi}\nabla v_\epsi\cdot\nabla \varphi+p\chi\intoTomega\frac{u_\epsi^p}{(1+\epsi u_\epsi)v_\epsi}\nabla v_\epsi\cdot\nabla\varphi-2(1-p)\chi\intoTomega\frac{u_\epsi^p}{v_\epsi}\nabla v_\epsi\cdot\nabla\varphi.\nonumber
\end{align}
\end{lemma}

\begin{bew}
We start using the first equation of \eqref{KSeps} and multiple integrations by parts to compute
\begin{align}\label{eq:testing-1}
\intomega\partial_t(u_\epsi^p\varphi)-\intomega u_\epsi^p\varphi_t&=p\intomega u_\epsi^{p-1}\Big[\Delta u_\epsi-\chi\nabla\cdot\Big(\frac{u_\epsi}{(1+\epsi u_\epsi)v_\epsi}\nabla v_\epsi\Big)\Big]\varphi\nonumber\\
&=-p\intomega \nabla\big(u_\epsi^{p-1}\varphi\big)\cdot\Big[\nabla u_\epsi-\chi\Big(\frac{u_\epsi}{(1+\epsi u_\epsi)v_\epsi}\nabla v_\epsi\Big)\Big]\nonumber\\
&=\frac{4(1-p)}{p}\intomega |\nabla u^\frac{p}{2}_\epsi|^2\varphi+\intomega u_\epsi^{p} \Delta\varphi\\
&\qquad-(1-p)\chi\intomega\nabla\Phi_{\epsi}(u_\epsi)\cdot\frac{\nabla v_\epsi}{v_\epsi}\varphi+p\chi\intomega\frac{u_\epsi^{p}}{(1+\epsi u_\epsi)v_\epsi}\nabla v_\epsi\cdot\nabla\varphi\quad\text{for all }t>0,\nonumber
\end{align}
where we rewrote $u_\epsi^{p-2}|\nabla u_\epsi|^2=\frac{4}{p^2}|\nabla u_\epsi^{\frac{p}{2}}|^2$. Integrating the integral containing $\nabla\Phi_{\epsi}(u_\epsi)$ once more by parts, relying on the second equation of \eqref{KSeps} to express $\Delta v_\epsi=v_\epsi-u_\epsi$ and making use of the fact that for all $\epsi\in(0,1)$ we have $\Phi_\epsi(u_\epsi)\leq u_\epsi^p$ for all $t>0$, we see that
\begin{align}\label{eq:testing-2}
\intomega\nabla\Phi_{\epsi}(u_\epsi)\cdot\frac{\nabla v_\epsi}{v_\epsi}\varphi
&=-\intomega\Phi_\epsi(u_\epsi)\Big[\frac{v_\epsi-u_\epsi}{v_\epsi}\Big]\varphi+\intomega\Phi_\epsi(u_\epsi)\frac{|\nabla v_\epsi|^2}{v_\epsi^2}\varphi-\intomega\frac{\Phi_\epsi(u_\epsi)}{v_\epsi}\nabla v_\epsi\cdot\nabla \varphi\nonumber\\
&\leq-\intomega\Phi_\epsi(u_\epsi)\varphi+\intomega \frac{u_\epsi^{p+1}}{v_\epsi}\varphi+\intomega u_\epsi^p\frac{|\nabla v_\epsi|^2}{v_\epsi^2}\varphi-\intomega\frac{\Phi_\epsi(u_\epsi)}{v_\epsi}\nabla v_\epsi\cdot\nabla \varphi
\end{align}
holds for all $t>0$. Now, to get rid of the term quadratic in $\frac{\nabla v_\epsi}{v_\epsi}$, we test the second equation in \eqref{KSeps} with $\frac{u_\epsi^p}{v_\epsi}\varphi$ -- which again due to Lemma \ref{lem:lowerbound-veps} is an admissible test function -- and integrate by parts to obtain
\begin{align*}
\intomega u_\epsi^p\frac{|\nabla v_\epsi|^2}{v_\epsi^2}\varphi-\intomega(\nabla u_\epsi^p\cdot \nabla v_\epsi)\frac{\varphi}{v_\epsi}-\intomega \frac{u_\epsi^p}{v_\epsi}\nabla v_\epsi\cdot\nabla \varphi-\intomega u_\epsi^p\varphi+\intomega\frac{u_\epsi^{p+1}}{v_\epsi}\varphi=0
\end{align*}
for all $t>0$, which by rewriting $\nabla u_\epsi^p=2 u_\epsi^{\frac{p}{2}}\nabla u_\epsi^{\frac{p}{2}}$ implies
\begin{align*}
\intomega u_\epsi^p\frac{|\nabla v_\epsi|^2}{v_\epsi^2}\varphi
&=\intomega u_\epsi^p\varphi-\intomega \frac{u_\epsi^{p+1}}{v_\epsi}\varphi+2\intomega \frac{u_\epsi^\frac{p}{2}}{v_\epsi}(\nabla u_\epsi^\frac{p}{2}\cdot\nabla v_\epsi)\varphi+\intomega \frac{u_\epsi^p}{v_\epsi}\nabla v_\epsi\cdot\nabla \varphi\quad\text{for all }t>0.
\end{align*} 
Decomposing the integral containing the mixed derivatives in the same manner as in the proof of Lemma \ref{lem:compability-class-sol}, we see that this readily implies
\begin{align}\label{eq:testing-3}
-(1-p)\chi\intomega u_\epsi^p\frac{|\nabla v_\epsi|^2}{v_\epsi^2}\varphi=&-2(1-p)\chi\intomega u_\epsi^p\varphi+2(1-p)\chi\intomega\frac{u_\epsi^{p+1}}{v_\epsi}\varphi+4(1-p)\chi\intomega\Big|\nabla u_\epsi^{\frac{p}{2}}-\frac{u_\epsi^{\frac{p}{2}}}{2v_\epsi}\nabla v_\epsi\Big|^2\varphi\nonumber\\
&-4(1-p)\chi\intomega|\nabla u_\epsi^{\frac{p}{2}}|^2\varphi-2(1-p)\chi\intomega\frac{u_\epsi^p}{v_\epsi}\nabla v_\epsi\cdot\nabla \varphi
\end{align}
for all $t>0$. Finally, a combination of \eqref{eq:testing-1}--\eqref{eq:testing-3} completes the proof upon integration over $(0,T)$.
\end{bew}

As an immediate consequence of the differential inequality provided by the preceding lemma we obtain the following spatio-temporal estimates for suitable values of $p\in(0,1)$.
\begin{lemma}\label{lem:st-bounds}
For $\chi>0$ let $p\in\big(0,1)$ satisfy $\chi<\frac{1}p$. Then for each $T>0$ there exists $C(p,T)>0$ such that
\begin{align}\label{eq:st-bound-nabla_ueps_p2}
\intoTomega|\nabla u_\epsi^\frac{p}{2}|^2\leq C(p,T)
\end{align}
and
\begin{align}\label{eq:st-bound-squares}
\intoTomega\Big|\nabla u_\epsi^\frac{p}{2}-\frac{u_\epsi^{\frac{p}{2}}}{2v_\epsi}\nabla v_\epsi\Big|^2\leq C(p,T)
\end{align}
as well as
\begin{align}\label{eq:st-bound_ueps_p+1_veps}
\intoTomega\frac{u_\epsi^{p+1}}{v_\epsi}\leq C(p,T)
\end{align}
and
\begin{align}\label{eq:st-bound_ueps_p_nabla_veps_l2}
\intoTomega \frac{u_\epsi^p}{v_\epsi^2}|\nabla v_\epsi|^2\leq C(p,T)
\end{align}
for all $\epsi\in(0,1)$.
\end{lemma}

\begin{bew}
By an application of Lemma \ref{lem:testing-equation} to $\varphi\equiv1$, by the positivity of $\Phi_\epsi(u_\epsi)$ we see that
\begin{align}\label{eq:st-bounds-eq1}
\frac{4(1-p)(1-p\chi)}{p}\intoTomega|\nabla u_\epsi^{\frac{p}{2}}|^2&+4(1-p)\chi\intoTomega\Big|\nabla u_\epsi^\frac{p}{2}-\frac{u_\epsi^\frac{p}{2}}{2v_\epsi}\nabla v_\epsi\Big|^2\\
&-2(1-p)\chi\intoTomega u_\epsi^p+(1-p)\chi\intoTomega\frac{u_\epsi^{p+1}}{v_\epsi}\leq \intomega u_\epsi^p(\cdot,T)-\intomega u_0^p\nonumber
\end{align}
holds for all $\epsi\in(0,1)$. Due to $0<\chi$ and $p<\min\{1,\frac{1}{\chi}\}$, this immediately implies \eqref{eq:st-bound-nabla_ueps_p2}--\eqref{eq:st-bound_ueps_p+1_veps} in view of Lemma \ref{lem:l1-bounds}. To see that also \eqref{eq:st-bound_ueps_p_nabla_veps_l2} holds, we make use of the strict positivity of $v_\epsi$ ensured in Lemma \ref{lem:lowerbound-veps} and test the second equation with $\frac{u_\epsi^{p}}{v_\epsi}$ to obtain that
\begin{align*}
0=\intomega\Delta v_\epsi\frac{u_\epsi^p}{v_\epsi}-\intomega u_\epsi^p+\intomega\frac{u_\epsi^{p+1}}{v_\epsi}=\intomega u_\epsi^p\frac{|\nabla v_\epsi|^2}{v_\epsi^2}-p\intomega\frac{u_\epsi^{p-1}}{v_\epsi}\nabla u_\epsi\cdot\nabla v_\epsi-\intomega u_\epsi^p+\intomega\frac{u_\epsi^{p+1}}{v_\epsi}
\end{align*}
holds for all $t>0$ and all $\epsi\in(0,1)$. Thus, integrating over $(0,T)$, rewriting $p u_\epsi^{p-1}\nabla u_\epsi=2u_\epsi^\frac{p}{2}\nabla u_\epsi^{\frac{p}{2}}$ and applying Young's inequality shows that
\begin{align*}
\intoTomega u_\epsi^p\frac{|\nabla v_\epsi|^2}{v_\epsi^2}\leq 4\intoTomega |\nabla u_\epsi^\frac{p}{2}|^2-2\intoTomega\frac{u_\epsi^{p+1}}{v_\epsi}+2\intoTomega u_\epsi^p\quad\text{for all }\epsi\in(0,1),
\end{align*}
which readily implies \eqref{eq:st-bound_ueps_p_nabla_veps_l2} in view of \eqref{eq:st-bound-nabla_ueps_p2} and $p<1$ combined with Lemma \ref{lem:l1-bounds}.
\end{bew}

The boundedness information on $\intoTomega\frac{u_\epsi^{p+1}}{v_\epsi}$ obtained in the previous lemma is the crucial ingredient in improving the regularity of $u_\epsi$. Here $p\in(0,1)$ must not be too small leading to the main reason for the restriction $\chi<\frac{n}{n-2}$.
\vspace*{-2pt}
\begin{lemma}\label{lem:st-bound-ueps_r}
For $0<\chi<\frac{n}{n-2}$ let $p\in(0,1)$ satisfy $\chi<\frac{1}p<\frac{n}{n-2}$. Then there exists some $r>1$ such that for any $T>0$ there exists $C(T)>0$ such that
\begin{align*}
\intoTomega u_\epsi^{r}\leq C(T)\quad\text{for all }\epsi\in(0,1).
\end{align*}
\end{lemma}

\begin{bew}
Given $0<\chi<\frac{n}{n-2}$ and $p\in(0,1)$ such that $\chi<\frac{1}{p}<\frac{n}{n-2}$ we can fix $r\in(1,p+1)$ satisfying $r<\frac{n(p+1)}{2n-2}$ and make use of Young's inequality to find $C_1>0$ such that
\vspace*{-2.5pt}
\begin{align*}
\intoTomega u_\epsi^{r}=\intoTomega\Big(\frac{u_\epsi^{p+1}}{v_\epsi}\Big)^\frac{r}{p+1}v_\epsi^{\frac{r}{p+1}}\leq \intoTomega\frac{u_\epsi^{p+1}}{v_\epsi}+C_1\intoTomega v_\epsi^{\frac{r}{p+1-r}}
\end{align*}
holds for all $\epsi\in(0,1)$. Since $r<\frac{n(p+1)}{2n-2}$ implies $\frac{r}{p+1-r}<\frac{n}{n-2}$, we can make use of Lemma \ref{lem:veps-bounds} and \eqref{eq:st-bound_ueps_p+1_veps} to find $C_2>0$ and $C_3>0$ satisfying
\vspace*{-2.5pt}
\begin{align*}
\intoTomega u_\epsi^r\leq C_3+C_1C_2 T\quad\text{for all }\epsi\in(0,1).\tag*{\qedhere}
\end{align*}
\end{bew}

\vspace*{-12pt}
\subsection{Time regularity of \texorpdfstring{$u_\epsi$}{u}}\label{sec4-2:timereg}

In pursuance of convergence properties suitable for our definition of generalized solutions we will rely on an Aubin-Lions type lemma for which we will require some additional information on the time regularity of our approximate solutions. We will therefore make use of some of the previously established a priori estimates to supplement our current repertoire of estimates with the following lemma.
\begin{lemma}\label{lem:timereg-ueps}
Assume $\chi>0$ and let $p\in(0,1)$ satisfy $\chi<\frac{1}p$. Then for all $T>0$ there exists $C(T)>0$ such that
\vspace*{-1pt}
\begin{align*}
\intoT\Big\|\partial_t\big(u_\epsi(\cdot,t)+1\big)^\frac{p}{2}\Big\|_{(W_0^{1,\infty}\!\left(\Omega\right))^*}\intd t\leq C(T)\quad\text{for all }\epsi\in(0,1).
\end{align*}
\end{lemma}

\begin{bew}
For fixed $\psi\in C_0^\infty\!\left(\Omega\right)$ such that $\|\psi\|_{\W[1,\infty]}\leq1$ we make use of the first equation in \eqref{KSeps} and integration by parts to obtain
\begin{align*}
\bigg|\intomega\partial_t(u_\epsi+1)^\frac{p}{2}\psi\bigg|&=\bigg|\frac{p}{2}\intomega(u_\epsi+1)^\frac{p-2}{2}\Big[\Delta u_\epsi-\chi\nabla\cdot\Big(\frac{u_\epsi}{(1+\epsi u_\epsi)v_\epsi}\nabla v_\epsi\Big)\Big]\psi\bigg|\\
&=\bigg|\frac{p(2-p)}{4}\intomega(u_\epsi+1)^\frac{p-4}{2}|\nabla u_\epsi|^2\psi-\frac{p}{2}\intomega(u_\epsi+1)^\frac{p-2}{2}\nabla u_\epsi\cdot\nabla\psi\\
&\qquad-\frac{p(2-p)\chi}{4}\intomega\frac{(u_\epsi+1)^\frac{p-4}{2}u_\epsi}{(1+\epsi u_\epsi)v_\epsi}(\nabla u_\epsi\cdot\nabla v_\epsi)\psi+\frac{p\chi}{2}\intomega\frac{(u_\epsi+1)^\frac{p-2}{2}u_\epsi}{(1+\epsi u_\epsi)v_\epsi}\nabla v_\epsi\cdot\nabla\psi\bigg|
\end{align*}
for all $t>0$ and $\epsi\in(0,1)$.
Having in mind the obvious estimates $u_\epsi\leq(u_\epsi+1)$, $(u_\epsi+1)^{-\frac{p}{2}}\leq 1$, and $(1+\epsi u_\epsi)^{-1}\leq1$ in $\Omega\times(0,T)$, we can draw on the fact that $\|\psi\|_{\W[1,\infty]}\leq1$ and Young's inequality to see that
\begin{align}\label{eq:time-reg-est}
\bigg|\intomega\partial_t(u_\epsi+1)^\frac{p}{2}\psi\bigg|&\leq\frac{p(2-p)}{4}\intomega(u_\epsi+1)^{p-2}|\nabla u_\epsi|^2+\frac{p}{4}\intomega(u_\epsi+1)^{p-2}|\nabla u_\epsi|^2+\frac{p|\Omega|}{4}\nonumber\\
&\qquad+\frac{p(2-p)\chi}{8}\intomega (u_\epsi+1)^{p-2}|\nabla u_\epsi|^2+\frac{p(2-p)\chi}{8}\intomega\frac{|\nabla v_\epsi|^2}{v_\epsi^2}\\&\qquad\qquad+\frac{p\chi}{4}\intomega(u_\epsi+1)^p+\frac{p\chi}{4}\intomega\frac{|\nabla v_\epsi|^2}{v_\epsi^2}\nonumber
\end{align}
holds for all $t>0$ and $\epsi\in(0,1)$. Since $p<1$, for all $t>0$ and $\epsi\in(0,1)$ we have that
\begin{align*}
\intomega(u_\epsi+1)^{p-2}|\nabla u_\epsi|^2\leq \intomega u_\epsi^{p-2}|\nabla u_\epsi|^2=\frac{4}{p^2}\intomega|\nabla u_\epsi^\frac{p}{2}|^2
\quad\text{ and }\quad 
\intomega(u_\epsi+1)^p\leq\intomega u_0+|\Omega|
\end{align*}
in view of Lemma \ref{lem:l1-bounds}, whereas Lemma \ref{lem:ln-nab-veps_l2} shows that
\begin{align*}
\intomega\frac{|\nabla v_\epsi|^2}{v_\epsi^2}\leq|\Omega|\quad\text{for all }t>0\text{ and }\epsi\in(0,1).
\end{align*}
Thus, a combination of these three estimates with \eqref{eq:time-reg-est} provides $C_1>0$ such that for all $\epsi\in(0,1)$,
\begin{align*}
\Big\|\partial_t\big(u_\epsi(\cdot,t)+1\big)^\frac{p}{2}\Big\|_{(W_0^{1,\infty}\!\left(\Omega\right))^*}\leq C_1\intomega|\nabla u_\epsi(\cdot,t)^\frac{p}{2}|^2+C_1\quad\text{for all }t\in(0,T),
\end{align*}
which in conjunction with \eqref{eq:st-bound-nabla_ueps_p2} of Lemma \ref{lem:st-bounds} completes the proof upon an integration over $t\in(0,T)$.
\end{bew}

\subsection{Convergence properties}\label{sec4-3:conv}

From the above estimates we can now extract a subsequence along which we may pass to the limit in a way suitable for our setting.   
\begin{lemma}\label{lem:convergence}
Assume $0<\chi<\frac{n}{n-2}$ and let $p\in(0,1)$ satisfy $\chi<\frac{1}p<\frac{n}{n-2}$. Then there exist $(\epsi_j)_{j\in\N}$ and functions $u$ and $v$ defined on $\Omega\times(0,\infty)$ such that $\epsi_j\searrow0$ as $j\to\infty$, that $u\geq0$ and $v\geq0$ a.e. in $\Omega\times(0,\infty)$, and that
\begin{align}
u_\epsi&\to u&&\text{in }\LSploc{1}{\bomega\times[0,\infty)}&&\text{and a.e. in }\Omega\times(0,\infty),\label{eq:conv_ueps_l1}\\
\nabla u_\epsi^\frac{p}{2}&\wto \nabla u^\frac{p}{2}&&\text{in }\LSploc{2}{\bomega\times[0,\infty)},&&\label{eq:conv_nab-ueps_wl2}\\
v_\epsi&\to v&&\text{in }\LSploc{1}{[0,\infty);\W[1,1]}&&\text{and a.e. in }\Omega\times(0,\infty),\label{eq:conv_veps_l1}\\
\frac{u_\epsi^{p+1}}{v_\epsi}&\to\frac{u^{p+1}}{v}&&\text{in }\LSploc{1}{\bomega\times[0,\infty)},&&\label{eq:conv_ueps-p+1-veps_l1}\\
u_\epsi^\frac{p}{2}\frac{\nabla v_\epsi}{v_\epsi}&\wto u^\frac{p}{2}\frac{\nabla v}{v}&&\text{in }\LSploc{2}{\Omega\times(0,\infty)},&&\label{eq:conv_ueps-nab-lnveps_wl2}
\end{align}
as $\epsi=\epsi_j\searrow0$. Moreover,
\begin{align}\label{eq:mass_cons_u}
\intomega u(\cdot,t)=\intomega u_0\qquad\text{for a.e. }t>0.
\end{align}
\end{lemma}

\begin{bew}
Intending to employ an Aubin-Lions type argument to obtain a first convergence information for $u_\epsi$, we fix any $p\in(0,1)$ such that $\chi<\frac{1}{p}<\frac{n}{n-2}$ and combine Lemma \ref{lem:st-bounds} with Lemma \ref{lem:l1-bounds} and Lemma \ref{lem:timereg-ueps} to find that
\begin{align*}
\Big((u_\epsi+1)^\frac{p}{2}\Big)_{\epsi\in(0,1)}\ \text{ is bounded in }\LSploc{2}{[0,\infty);\W[1,2]}
\end{align*}
and that
\begin{align*}
\Big(\partial_t(u_\epsi+1)^\frac{p}{2}\Big)_{\epsi\in(0,1)}\ \text{ is bounded in }L_{loc}^1\big([0,\infty); (W_0^{1,\infty}\!\left(\Omega\right))^*\big).
\end{align*}
Hence, we can invoke an Aubin-Lions lemma (\cite[Corollary 8.4]{Sim87}) to infer the existence of $(\epsi_j)_{j\in\N}\subset(0,1)$ such that $\epsi_j\searrow0$ as $j\to\infty$, that
\begin{align}\label{eq:conveq1}
u_\epsi^\frac{p}{2}\to u^\frac{p}{2}\quad\text{in }\LSploc{2}{\bomega\times[0,\infty)}\text{ and a.e. in }\Omega\times(0,\infty)\text{ as }\epsi=\epsi_j\searrow0
\end{align}
and such that \eqref{eq:conv_nab-ueps_wl2} holds with some nonnegative function $u$ defined on $\Omega\times(0,\infty)$. Now, relying on the a.e. convergence of $u_\epsi^\frac{p}{2}$ and the equi-integrability property of $\{u_{\epsi_j}^{1+s}\}_{j\in\N}$ for some small $s>0$ contained in Lemma \ref{lem:st-bound-ueps_r} we may employ the Vitali convergence theorem to find 
\begin{align}\label{eq:convergence-eq1}
u_\epsi\to u\quad\text{in }\LSploc{1+s}{\bomega\times[0,\infty)}\ \text{as }\epsi=\epsi_j\searrow0\text{ for some small }s>0,
\end{align}
implying that also \eqref{eq:conv_ueps_l1} holds, whereupon \eqref{eq:mass_cons_u} follows from Lemma \ref{lem:l1-bounds}. Since from standard elliptic theory (e.g. \cite[Theorem 9.32]{Brezis_funcana}) we know that for all $\epsi\in(0,1)$ and all $r>11$ we have $\|v_\epsi\|_{\W[2,r]}\leq C\|u_\epsi\|_{\Lo[r]}$ with $C>0$, \eqref{eq:convergence-eq1} readily implies that there exists some nonnegative $v$ defined on $\Omega\times(0,\infty)$ such that \eqref{eq:conv_veps_l1} holds. For proving \eqref{eq:conv_ueps-p+1-veps_l1} we pick $q>p$ such that still $\chi<\frac{1}{q}$ holds and see that by Lemma \ref{lem:st-bounds} there exists $C_2>0$ such that with $K_1>0$ taken from  Lemma \ref{lem:lowerbound-veps} we have
\begin{align*}
\intoTomega \Big(\frac{u_\epsi^{p+1}}{v_\epsi}\Big)^{\frac{q+1}{p+1}}=\intoTomega\frac{u_\epsi^{q+1}}{v_\epsi}\frac{1}{v_\epsi^{\frac{q-p}{p+1}}}\leq K_1^\frac{p-q}{p+1}\intoTomega\frac{u_\epsi^{q+1}}{v_\epsi} \leq K_1^{\frac{p-q}{p+1}}C_2
\end{align*}
for all $\epsi\in(0,1)$. Hence, an application of the Vitali convergence theorem proves \eqref{eq:conv_ueps-p+1-veps_l1}. To verify \eqref{eq:conv_ueps-nab-lnveps_wl2}, we first note that in view of \eqref{eq:st-bound_ueps_p_nabla_veps_l2} there exists some $w\in\LSp{2}{\Omega\times(0,\infty)}$ such that (upon choice of a suitable subsequence)
\begin{align}\label{eq:conv_ueps-nab-lnveps-to-some-w}
u_\epsi^\frac{p}{2}\frac{\nabla v_\epsi}{v_\epsi}\wto w\quad\text{in }\LSploc{2}{\Omega\times(0,\infty)}\ \text{as }\epsi=\epsi_j\searrow0.
\end{align}
Furthermore, due to \eqref{eq:conveq1} we have $u_\epsi^\frac{p}{2}\to u^\frac{p}{2}$ in $\LSploc{2}{\Omega\times(0,\infty)}$ and $\frac{\nabla v_\epsi}{v_\epsi}\wto\frac{\nabla v}{v}$ in $\LSploc{2}{\Omega\times(0,\infty)}$ in light of the precompactness property implied by Lemma \ref{lem:ln-nab-veps_l2}, and thus
\begin{align*}
u_\epsi^\frac{p}{2}\frac{\nabla v_\epsi}{v_\epsi}\wto u^\frac{p}{2}\frac{\nabla v}{v}\quad\text{in }\LSploc{1}{\Omega\times(0,\infty)}\ \text{as }\epsi=\epsi_j\searrow0,
\end{align*}
which ensures that $w=u^\frac{p}{2}\frac{\nabla v}{v}$ a.e. in $\Omega\times(0,T)$, and due to \eqref{eq:conv_ueps-nab-lnveps-to-some-w} hence shows \eqref{eq:conv_ueps-nab-lnveps_wl2}.
\end{bew}

\setcounter{equation}{0} 
\section{Proof of Theorem \ref{theo:globgensol}}\label{sec5:prooftheo}

In order to verify the crucial positivity properties demanded in the Definition \ref{def:glob-gen-sol} we want to find some lower bound for $\intomega \ln u_\epsi$. To this end we will state two technical lemmas which have been proven in \cite{LanWin2017} and prepare a comparison argument for a differential inequality.  
\begin{lemma}\label{lem:diffineq-lem}
Let $a>0$, $b>0$, $T>0$ and let $y:(0,T)\to\R$ be a continuously differentiable function satisfying
\begin{align*}
y'(t)\leq -ay^2(t)+b\qquad\text{for all }t\in(0,T),\text{ for which }y(t)>0.
\end{align*} 
Then
\begin{align*}
y(t)\leq\sqrt{\frac{b}{a}}\coth\big(\sqrt{ab}t\big)\qquad\text{for all }t\in(0,T).
\end{align*}
\end{lemma}
\begin{bew}
We refer the reader to \cite[Lemma 8.3]{LanWin2017} for the proof.
\end{bew}

In addition to the previous comparison lemma we will also make use of the following auxiliary lemma which was given in \cite[Lemma 8.4]{LanWin2017} -- generalizing a result proven in \cite[Lemma 4.3]{TaoWin_JDE15} to non-convex domains.
\begin{lemma}\label{lem:lanwinpart2}
Let $\eta>0$. Then there exists $C>0$ such that every positive function $\varphi\in \CSp{1}{\bomega}$ fulfilling
\begin{align*}
|\{x\in\Omega:\,\varphi(x)>\delta\}|>\eta
\end{align*}
for some $\delta>0$ satisfies
\begin{align*}
\intomega\frac{|\nabla \varphi|^2}{\varphi^2}>C\left(\intomega\ln\frac{\delta}{\varphi}\right)^2\quad\text{or}\quad\intomega\ln\frac{\delta}{\varphi}<0.
\end{align*}
\end{lemma}
\begin{bew}
This is \cite[Lemma 8.4]{LanWin2017}.
\end{bew}

Relying on the previous two lemmata we can now build on the ideas from \cite[Lemma 8.5]{LanWin2017} obtain the following.
\begin{lemma}\label{lem:ln-ueps_lower_bound}
There exists $T>0$ such that for every $t\in(0,T)$ the inequality 
\begin{align*}
\inf_{\epsi\in(0,1)}\intomega\ln u_\epsi(\cdot,t)>-\infty
\end{align*}
is valid.
\end{lemma}

\begin{bew}
For $t\in(0,\infty)$ and $\epsi\in(0,1)$ we let $M_\epsi(t):=\sup_{s\in[0,t]}\|u_\epsi(\cdot,s)\|_{\Lo[\infty]}$. Picking $q\in(n,\infty)$, we obtain from standard elliptic regularity theory that $\|v_\epsi\|_{\W[2,q]}\leq C_1\|u_\epsi\|_{\Lo[q]}$ for all $t\in(0,\infty)$ with some $C_1>0$ and hence, by the Sobolev embedding theorem, that
\begin{align*}
\|\nabla v_\epsi(\cdot,t)\|_{\Lo[\infty]}\leq C_2 M_\epsi(t) \quad\text{for all }t\in(0,\infty).
\end{align*}
Invoking the well-known smoothing estimates for the Neumann heat semigroup (e.g. \cite[Lemma 2.1]{caolan16_smalldatasol3dnavstokes} or \cite[Lemma 1.3]{win10jde}) we find $C_3>0$ such that
\begin{align*}
\|u_\epsi(\cdot,t)\|_{\Lo[\infty]}&\leq \|u_0\|_{\Lo[\infty]}+C_3\chi\int_0^t\Big(1+(t-s)^{-\frac{1}{2}-\frac{n}{2q}}\Big)\Big\|\frac{u_\epsi(\cdot,s)}{\big(1+\epsi u_\epsi(\cdot,s)\big)v_\epsi(\cdot,s)}\nabla v_\epsi(\cdot,s)\Big\|_{\Lo[q]}\intd s\\
&\leq \|u_0\|_{\Lo[\infty]}+\frac{C_3\chi}{K_1}C_2 M_\epsi(t)\int_0^t\Big(1+(t-s)^{-\frac{1}{2}-\frac{n}{2q}}\Big)\|u_\epsi(\cdot,s)\|_{\Lo[q]}\intd s\\
&\leq \|u_0\|_{\Lo[\infty]}+\frac{C_3\chi}{K_1}C_2 M_\epsi(t)\|u_0\|_{\Lo[1]}^\frac{1}{q} M_\epsi(t)^\frac{q-1}{q}\int_0^t\Big(1+(t-s)^{-\frac{1}{2}-\frac{n}{2q}}\Big)\intd s
\end{align*}
for all $t\in(0,\infty)$ and all $\epsi\in(0,1)$, with $K_1>0$ given by Lemma \ref{lem:lowerbound-veps}. From this we infer the existence of $C_4>0$ such that
\begin{align*}
M_\epsi(t)\leq \|u_0\|_{\Lo[\infty]}+C_4 M_\epsi(t) M_\epsi^\frac{q-1}{q}(t)\big(t+t^{\frac{1}{2}-\frac{n}{2q}}\big)\quad\text{for all }t\in(0,\infty),
\end{align*}
which, in view of the fact that for all $a,b,\in[0,\infty)$ and $\gamma\in(0,1)$ the inequality
\begin{align*}
\sup\{x\in[0,\infty)\,|\, x\leq a+bx^\gamma\}\leq\frac{a}{1-\gamma}+b^\frac{1}{1-\gamma}
\end{align*}
holds true, implies that
\begin{align*}
M_\epsi(t)\leq q\|u_0\|_{\Lo[\infty]}+\Big(C_4M_\epsi(t)\big(t+t^{\frac{1}{2}-\frac{n}{2q}}\big)\Big)^q\quad\text{for all }t\in(0,\infty).
\end{align*}
Letting $T_\epsi:=\sup\left\{t\in(0,\infty)\,|\, M_\epsi(t)\leq q\|u_0\|_{\Lo[\infty]}+1\right\}$ we see that certainly
\begin{align*}
T_\epsi>T:=\min\left\{1,\frac{1}{2C_4\cdot 2(q\|u_0\|_{\Lo[\infty]}+1)}\right\},
\end{align*}
so that for all $\epsi\in(0,1)$ we may estimate 
\begin{align*}
\|u_\epsi(\cdot,t)\|_{\Lo[\infty]}\leq q\|u_0\|_{\Lo[\infty]}+1=:M\quad\text{for all }t\in(0,T).
\end{align*}
Now, since for $\delta:=\frac{1}{2|\Omega|}\intomega u_0$ and $\eta:=\frac{1}{2M}\intomega u_0$ we have 
\begin{align*}
\intomega u_0=\intomega u_\epsi(\cdot,t)=\int_{\{u_\epsi(\cdot,t)\geq\delta\}}u_\epsi(\cdot,t)+\int_{\{u_\epsi(\cdot,t)<\delta\}} u_\epsi(\cdot,t)\leq M|\{u_\epsi(\cdot,t)\geq\delta\}|+\delta|\Omega|,
\end{align*}
we conclude $|\{u_\epsi(\cdot,t)\geq\delta\}|\geq\eta$ for every $t\in(0,T)$ and each $\epsi\in(0,1)$, so that from Lemma \ref{lem:lanwinpart2} we obtain some $C_5>0$ such that a combination of Lemma \ref{lem:lanwinpart2} with Lemma \ref{lem:ln-nab-veps_l2} shows that
\begin{align*}
\frac{\intd}{\intd t}\Big(\intomega\ln\frac{\delta}{u_\epsi}\Big)&=-\intomega\frac{|\nabla u_\epsi|^2}{u_\epsi^2}+\chi\intomega\frac{1}{(1+\epsi u_\epsi)u_\epsi v_\epsi}\nabla u_\epsi\cdot\nabla v_\epsi\\
&\leq -\frac{1}{2}\intomega\frac{|\nabla u_\epsi|^2}{u_\epsi^2}+\frac{\chi^2}{2}\intomega\frac{|\nabla v_\epsi|^2}{v_\epsi^2}\leq -\frac{C_5}{2}\Big(\intomega\ln\frac{\delta}{u_\epsi}\Big)^2+\frac{\chi^2|\Omega|}{2}
\end{align*}
for every $t\in(0,T)$ where $\intomega\ln\frac{\delta}{u_\epsi(\cdot,t)}>0$. Therefore, an application of Lemma \ref{lem:diffineq-lem} completes the proof.
\end{bew}

The three preceding lemmas at hand we can now emulate the arguments featured in \cite[Lemma 8.6]{LanWin2017} to verify the essential positivity requirements appearing in Definition \ref{def:glob-gen-sol}. 
\begin{lemma}\label{lem:positivity_u_v}
Assume $0<\chi<\frac{n}{n-2}$ and let $p\in(0,1)$ be such that $\chi<\frac{1}p<\frac{n}{n-2}$. Then the functions $u$ and $v$ obtained in Lemma \ref{lem:convergence} satisfy $u>0$, and $v>0$ a.e. in $\Omega\times(0,\infty)$ as well as $u^\frac{p}{2}>0$ a.e. on $\romega\times(0,\infty)$.
\end{lemma}

\begin{bew}
The positivity of $v$ a.e. in $\Omega\times(0,\infty)$ follows from Lemma \ref{lem:lowerbound-veps} and \eqref{eq:conv_veps_l1}. For the positivity of $u$ a.e. in $\Omega\times(0,\infty)$ and a.e. on $\romega\times(0,\infty)$ we start by calculating
\begin{align}\label{eq:posi-eq1}
-\frac{\intd}{\intd t}\intomega \ln u_\epsi=-\intomega\frac{|\nabla u_\epsi|^2}{u_\epsi^2}+\chi\intomega\frac{1}{u_\epsi(1+\epsi u_\epsi)v_\epsi}\nabla u_\epsi\cdot\nabla v_\epsi\leq-\frac{1}{2}\intomega\frac{|\nabla u_\epsi|^2}{u_\epsi^2}+\frac{\chi^2}{2}\intomega\frac{|\nabla v_\epsi|^2}{v_\epsi^2}
\end{align}
for all $t>0$. Now, for fixed $\tau>0$ in view of Lemma \ref{lem:ln-ueps_lower_bound} we can find $\tau_0\in(0,\tau)$ such that
\begin{align}\label{eq:posi-eq2}
\inf_{\epsi\in(0,1)}\intomega\ln u_\epsi(\cdot,\tau_0)>-\infty,
\end{align}
which together with Lemma \ref{lem:ln-nab-veps_l2} shows upon integration of \eqref{eq:posi-eq1} that for any fixed $T>\tau$ we have
\begin{align}\label{eq:posi-eq3}
-\intomega\ln u_\epsi(\cdot,T)+\frac{1}{2}\int_{\tau_0}^t\!\intomega|\nabla \ln u_\epsi|^2\leq -\intomega \ln u_\epsi(\cdot,\tau_0)+(T-\tau_0)\frac{\chi^2|\Omega|}{2}
\end{align}
for all $t\in(\tau,T)$, with the right-hand side being bounded independently of $\epsi$ by virtue of \eqref{eq:posi-eq2}. Relying on the basic estimate $|\ln \xi|\leq 2\xi-\ln\xi$ for all $\xi>0$, we can first make use of Lemma \ref{lem:l1-bounds} and \eqref{eq:posi-eq3} to find $C_1>0$ such that $\intomega|\ln u_\epsi(\cdot,t)|\leq C_1$ for all $t\in(\tau,T)$. Afterwards, we invoke the Poincaré inequality to find $C_2>0$ such that $\|\varphi\|_{\W[1,2]}\leq C_2(\|\nabla \varphi\|_{\Lo[2]}+\|\varphi\|_{\Lo[1]})$ for all $\varphi\in\W[1,2]$ and conclude, again by \eqref{eq:posi-eq3}, that for every $\tau>0$ and $T>\tau$ there exists $C_3>0$ such that
\begin{align*}
\|\ln u_\epsi\|_{L^{2}\left((\tau,T);\W[1,2]\right)}\leq C_3\quad\text{for all }\epsi\in(0,1).
\end{align*}
In view of a weak compactness argument this means that we actually have $\ln u\in\LSploc{2}{(0,\infty);\W[1,2]}$, which readily entails $\ln u\in\LSploc{2}{\bomega\times(0,\infty)}$ and also $\ln u\in\LSploc{2}{\romega\times(0,\infty)}$ by a trace embedding theorem and thus proves the asserted positivity properties.
\end{bew}

Most of the requirements appearing in Definition \ref{def:glob-gen-sol} are prepared and all that is left is to combine the information presented in Lemma \ref{lem:diffineq-lem}, Lemma \ref{lem:convergence} and Lemma \ref{lem:positivity_u_v}.
\begin{proof}[\textbf{Proof of Theorem \ref{theo:globgensol}:}]
We fix $p\in(0,1)$ such that $\chi<\frac1p<\frac{n}{n-2}$, which in particular means that the requirements for Lemma \ref{lem:convergence} and Lemma \ref{lem:positivity_u_v} are satisfied. The regularity properties prescribed in \eqref{eq:sol_def_reg1} are satisfied according to \eqref{eq:conv_ueps_l1}, \eqref{eq:conv_veps_l1} and the regularity requirements featured in \eqref{eq:sol_def_reg2} are fulfilled in view of \eqref{eq:conv_nab-ueps_wl2}, \eqref{eq:conv_ueps-p+1-veps_l1}, and \eqref{eq:conv_ueps-nab-lnveps_wl2}. The positivity properties have been shown in Lemma \ref{lem:positivity_u_v}, whereas the mass identity \eqref{eq:mass} is valid due to \eqref{eq:mass_cons_u}. Since for all $\epsi\in(0,1)$ and arbitrary $\psi\in C_0^\infty\!\left(\bomega\times[0,\infty)\right)$ the global classical solution $(u_\epsi,v_\epsi)$ of \eqref{KSeps} satisfies
\begin{align*}
-\intinfomega\nabla v_\epsi\cdot\nabla \psi+\intinfomega v_\epsi\psi=\intinfomega u_\epsi\psi,
\end{align*}
we see that due to \eqref{eq:conv_ueps_l1} and \eqref{eq:conv_veps_l1} we may let $\epsi=\epsi_j\searrow0$ in each integral and obtain that \eqref{eq:sol_def_veq} holds and the only thing left is to verify \eqref{eq:sol_def_uineq}. To this end, we fix a nonnegative $\varphi\in C_0^\infty\!\left(\bomega\times[0,\infty)\right)$ satisfying $\frac{\partial\varphi}{\partial\nu}=0$ on $\romega\times(0,\infty)$ and $T>0$ such that $\varphi\equiv0$ in $\Omega\times[T,\infty)$. Invoking Lemma \ref{lem:testing-equation} shows that with $\Phi_\epsi$ as introduced in \eqref{eq:Phi-def},
\begin{align}\label{eq:proof-theo-1}
-&\intoTomega u_\epsi^p\varphi_t-\intomega u_0^p\varphi(\cdot,0)\nonumber\\
&\quad\geq \frac{4(1-p)(1-p\chi)}{p}\intoTomega|\nabla u_\epsi^\frac{p}{2}|^2\varphi+\intoTomega u_\epsi^p\Delta\varphi+4(1-p)\chi\intoTomega\Big|\nabla u_\epsi^\frac{p}{2}-\frac{u_\epsi^\frac{p}{2}}{2 v_\epsi}\nabla v_\epsi\Big|^2\varphi\nonumber\\
&\quad\ +(1-p)\chi\intoTomega\Phi_\epsi(u_\epsi)\varphi-2(1-p)\chi\intoTomega u_\epsi^p\varphi+(1-p)\chi\intoTomega\frac{u_\epsi^{p+1}}{v_\epsi}\varphi\\
&\quad\ +(1-p)\chi\intoTomega\frac{\Phi_\epsi(u_\epsi)}{v_\epsi}\nabla v_\epsi\cdot\nabla \varphi+p\chi\intoTomega\frac{u_\epsi^p}{(1+\epsi u_\epsi)v_\epsi}\nabla v_\epsi\cdot\nabla\varphi-2(1-p)\chi\intoTomega\frac{u_\epsi^p}{v_\epsi}\nabla v_\epsi\cdot\nabla\varphi.\nonumber
\end{align}
Relying on the facts that $p<1$ and $\Phi_\epsi(u_\epsi)\leq u_\epsi^p$ for all $\epsi\in(0,1)$ we see that by \eqref{eq:conv_ueps_l1} we have
\begin{align*}
-\intoTomega u_\epsi^p\varphi_t&\to -\intoTomega u^p\varphi_t,&
\intoTomega u_\epsi^p\Delta\varphi&\to \intoTomega u^p\Delta\varphi,\\
(1-p)\chi\intoTomega\Phi_\epsi(u_\epsi)\varphi&\to (1-p)\chi\intoTomega u^p\varphi,\quad\text{and}&
-2(1-p)\chi\intoTomega u_\epsi^p\varphi&\to-2(1-p)\chi\intoTomega u^p\varphi
\end{align*}
as $\epsi=\epsi_j\searrow0$, whereas \eqref{eq:conv_ueps-p+1-veps_l1} shows that
\begin{align*}
(1-p)\chi\intoTomega\frac{u_\epsi^{p+1}}{v_\epsi}\varphi\to(1-p)\chi\intoTomega\frac{u^{p+1}}{v}\varphi
\end{align*}
as $\epsi=\epsi_j\searrow0$. For the three integrals containing the spatial derivative of $v_\epsi$ we note that $\frac{u_\epsi^p}{v_\epsi}\nabla v_\epsi=\Big(u_\epsi^\frac{p}{2}\Big)\cdot\Big(\frac{u_\epsi^\frac{p}{2}}{v_\epsi}\nabla v_\epsi\Big)$, which according to \eqref{eq:conv_ueps_l1} and \eqref{eq:conv_ueps-nab-lnveps_wl2} implies
\begin{align*}
-2(1-p)\chi\intoTomega\frac{u_\epsi^p}{v_\epsi}\nabla v_\epsi\cdot\nabla \varphi&\to-2(1-p)\chi\intoTomega\frac{u^p}{v}\nabla v\cdot\nabla \varphi
\intertext{and, relying on the Vitali convergence theorem once more, also}
(1-p)\chi\intoTomega\frac{\Phi_\epsi(u_\epsi)}{v_\epsi}\nabla v_\epsi\cdot\nabla \varphi&\to(1-p)\chi\intoTomega\frac{u^p}{v}\nabla v\cdot\nabla\varphi,
\intertext{as well as}
p\chi\intoTomega\frac{u_\epsi^p}{(1+\epsi u_\epsi)v_\epsi}\nabla v_\epsi\cdot\nabla\varphi&\to p\chi\intoTomega\frac{u^p}{v}\nabla v\cdot\nabla\varphi,
\end{align*}
as $\epsi=\epsi_j\searrow0$. Finally, by the lower semicontinuity of the norm in $\LSp{2}{\Omega\times(0,T)}$ with respect to weak convergence it follows from \eqref{eq:conv_nab-ueps_wl2} and \eqref{eq:conv_ueps-nab-lnveps_wl2} that
\begin{align*}
\liminf_{\epsi_j\searrow0}\intoTomega|\nabla u_{\epsi_j}^\frac{p}{2}|^2\varphi&\geq\intoTomega|\nabla u^\frac{p}{2}|^2\varphi\\\text{and}\quad\liminf_{\epsi_j\searrow0}\intoTomega\Big|\nabla u_{\epsi_j}^\frac{p}{2}-\frac{u_{\epsi_j}^\frac{p}{2}}{2 v_{\epsi_j}}\nabla v_{\epsi_j}\Big|^2\varphi&\geq\intoTomega\Big|\nabla u^\frac{p}{2}-\frac{u^\frac{p}{2}}{2 v}\nabla v\Big|^2\varphi,
\end{align*}

so that passing to the limit each of the integrals in \eqref{eq:proof-theo-1} yields
\begin{align*}
-\intoTomega u^p&\varphi_t-\intomega u_0^p\varphi(\cdot,0)\nonumber\\
&\geq \frac{4(1-p)(1-p\chi)}{p}\intoTomega|\nabla u^\frac{p}{2}|^2\varphi+\intoTomega u^p\Delta\varphi+4(1-p)\chi\intoTomega\Big|\nabla u^\frac{p}{2}-\frac{u^\frac{p}{2}}{2v}\nabla v\Big|^2\varphi\nonumber\\
&\qquad -(1-p)\chi\intoTomega u^p\varphi+(1-p)\chi\intoTomega\frac{u^{p+1}}{v}\varphi-(1-2p)\chi\intoTomega\frac{u^p}{v}\nabla v\cdot\nabla\varphi,
\end{align*}
proving \eqref{eq:sol_def_uineq}, and thereby verifies that indeed $(u,v)$ is a global generalized solution in the sense of Definition \ref{def:glob-gen-sol}.
\end{proof}

\section*{Acknowledgements}
The author acknowledges support of the {\em Deutsche Forschungsgemeinschaft} in the context of the project
  {\em Analysis of chemotactic cross-diffusion in complex frameworks}. 


\footnotesize{

}
\end{document}